\documentclass[12pt]{article}

\usepackage[cp1251]{inputenc}
\usepackage[T2A]{fontenc}
\usepackage[english]{babel}
\usepackage{epsfig}
\usepackage{wrapfig}

\begin{document}

\centerline{\Large On discrete pseudohyperbolic attractors of Lorenz type}

\vspace{0.5cm}

\centerline{{\bf S.V. Gonchenko$^{1,2}$, A.S. Gonchenko$^{1}$,  A.O. Kazakov$^{2}$}}

\vspace{0.5cm}

{$^1$Lobachevsky State University of Nizhny Novgorod, Gagarina av. 23,  603950 Nizhny Novgorod, Russia} \\~\\
{$^2$National Research University Higher School of Economics, 25/12 Bolshaya Pecherskaya Ul., 603155 Nizhny Novgorod, Russia}


\vspace{0.5cm}

\centerline{\footnotesize{\em e-mail: sergey.gonchenko@mail.ru; agonchenko@mail.ru; kazakovdz@yandex.ru}}

\vspace{1.5cm}


\centerline{\textbf{Abstract.}}

\vspace{0.5cm}

{\footnotesize  We study geometrical and dynamical properties of the so-called discrete Lorenz-like attractors, that can be observed in three-dimensional diffeomorphisms. We propose new phenomenological scenarios of their appearance in one parameter families of such maps. We pay especially our attention to such a scenario that can lead to period-2 Lorenz-like attractors. These attractors have very interesting dynamical properties and we show that their crises can lead, in turn, to the emergence of pseudohyperbolic discrete Lorenz shape attractors of new types. We also show examples of all these attractors in three-dimensional generalized H\'enon maps. }


\section{Introduction}

In the theory of dynamical chaos, problems related to the study of strange attractors of multidimensional systems (with dimension of phase space $\geq 4$ for flows and $\geq 3$ for diffeomorphisms) are the most important.
Naturally, the solution of these problems
should be based on
fundamental results obtained in the theory of chaos in smaller dimensions.
First of all, it should be noted the results that are associated with the discovery of the Lorenz attractor \cite{Lor63}, the creation of its mathematical models
\cite{G76,ABS77,GW79,ABS82}, as well as with the proof of its robust chaoticity \cite{Sh81,SST93,Tucker99,OT17}. Together, these results, as well as studies on H\'enon-like attractors of two-dimensional diffeomorphisms \cite{H76,S79,CGP88,BenCar91,U95}, served as a kind of foundation on which the basic elements of the theory of discrete homoclinic attractors of three-dimensional maps were built \cite{GOST05,TS08,GGS12,GGKT14}.

Recall that a strange attractor of a diffeomorphism is called homoclinic if it contains only one saddle fixed point and, hence, entirely its unstable manifold.   In the case of two-dimensional maps, examples of discrete homoclinic attractors are well-known: the above H\'enon-like attractors, various Mira attractors \cite{Mira}, attractors in periodically perturbed systems with a homoclinic figure-8 of a saddle \cite{GSV13} etc.

%

In the present paper, we  consider three-dimensional diffeomorphisms which, on the one hand, are an independent subject of the theory of dynamical systems and, on the other hand, can be represented as Poincar\'e maps -- a quite effective tool for the study of four-dimensional flows.

As far as we know, the first explicit
examples of discrete homoclinic attractors in the case of three-dimensional diffeomorphisms were given in the paper \cite{GOST05}, in which they were found in three-dimensional maps of the form
\begin{equation}
\bar x = y,\;\; \bar y = z,\;\; \bar z = M_1 + B x + M_2 y - z^2,
\label{3dHM1}
\end{equation}
where $M_1,M_2$ and $B$ are parameters, $B$ is the Jacobian of map (\ref{3dHM1}). Note that map (\ref{3dHM1}) is a three-dimensional extension for the standard H\'enon map $\bar y = z,\; \bar z = M_1 + M_2 y - z^2$ (the latter is obtained at $B=0$), therefore, we will call map (\ref{3dHM1})  the {\em three-dimensional H\'enon map}.

\begin{figure}[ht]
\centerline{\epsfig{file=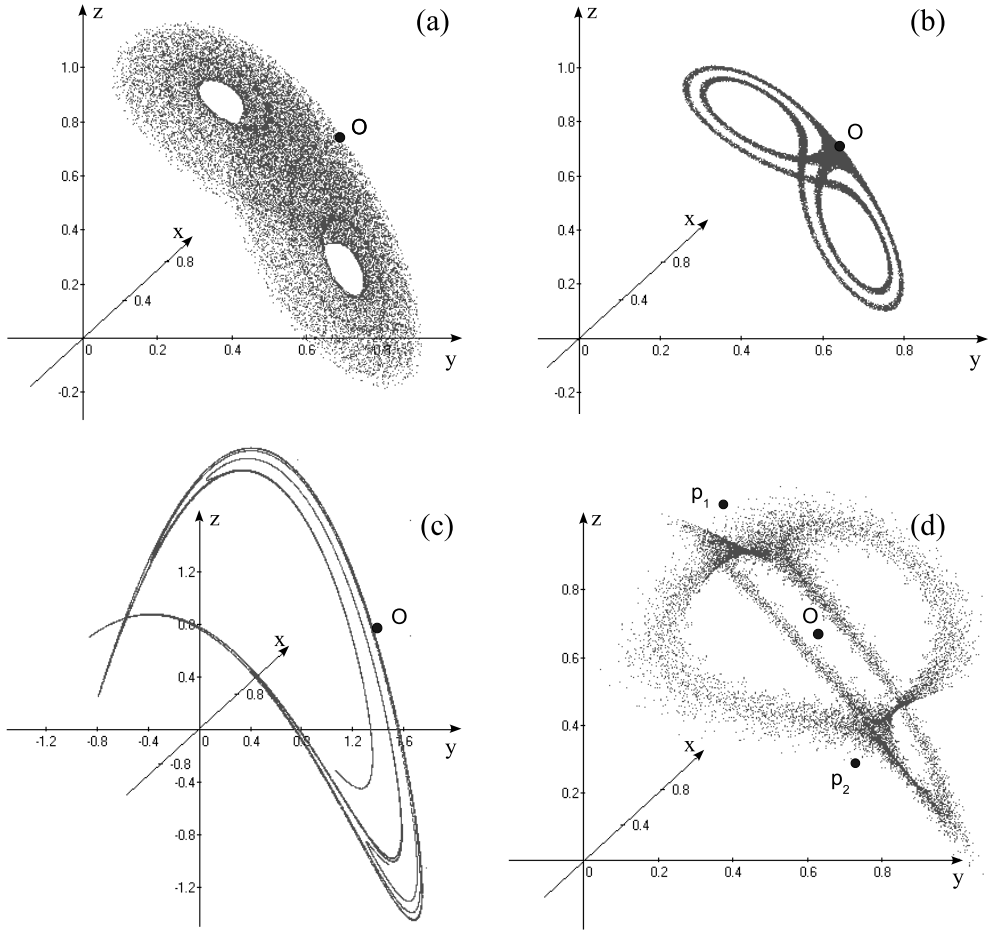, width=16cm
}}
\caption{{\footnotesize Plots of attractors of map (\ref{3dHM1}), based on the paper \cite{GOST05}. Figs.~(a)-(b): discrete Lorenz attractors at $B=0.7; M_1=0$ and (a) $M_2 = 0.85$ or (b) $M_2=0.815$. Fig.~(c): an attractor similar to the 2D H\'enon attractor, for $B=0.1; M_1=1.4;  M_2 = 0.2$. Fig.~(d): a homoclinic period-2 attractor (containing period-2 saddles $p_1$ and $p_2$ of type (2,1)), for $B= -0.95; M_1=1.77;  M_2 = - 0.925$.   }}
\label{fig:attrmap1}
\end{figure}

In Fig.~\ref{fig:attrmap1} phase portraits of some attractors of map (\ref{3dHM1}) are shown (here about $10^5$ forward iterations of a single point are taken for every plot).

The first two examples, shown in Figs.~\ref{fig:attrmap1}a,b, relate to
discrete homoclinic
attractors containing a saddle fixed point $O\left(x=y=z=B + M_2-1\right)$.
We note that the phase portraits of attractors in  Figs.~\ref{fig:attrmap1}a,b are very similar to the classical Lorenz attractors (without lacunae\footnote{Recall that in the case of the classical Lorenz attractor, a (trivial) lacuna is an open region (a hole inside the attractor) that contains a saddle limit cycle whose the stable manifold does not intersect the attractor \cite{ABS82}. Accordingly, in the case of discrete Lorenz attractor, a lacuna contains a saddle closed invariant curve with the same property. When varying parameters, lacunae  can appear and disappear due to the appearance/disappearance of homoclinic intersections. The corresponding bifurcations belong to the class of the so-called internal bifurcations of attractor. See more details in \cite{ABS82,Mal85,Mal03}.}
in Fig.~\ref{fig:attrmap1}a, and with a lacuna in Fig.~\ref{fig:attrmap1}b), therefore, these attractors were called in \cite{GGS12} {\em discrete Lorenz attractors}.

The attractor of Fig.~\ref{fig:attrmap1}c is also a discrete homoclinic attractor (it contains the fixed point $O(x=y=z\approx 0.88)$) which is similar to the two-dimensional H\'enon attractor \cite{Henon76}.
On the other hand, attractors of Figs.~\ref{fig:attrmap1}(a), (b) and (d) are essentially three-dimensional.
Besides, the attractor of
Fig.~\ref{fig:attrmap1}d
is nonorientable
(here $B= -0.95 <0$), and it is a homoclinic period-2 attractor, since it contains a period-2 saddle cycle $(p_1,p_2)$ (i.e. $p_2 = T(p_1)$, $p_1 = T(p_2)$). In the case under consideration, $p_1 = (x_0,y_0,x_0)$ and
$p_2=(y_0,x_0,y_0)$, where $x_0\approx 0.85, y_0 \approx 0.126$.
%

Note that map (\ref{3dHM1}) is a representative of a class of maps
\begin{equation}
\bar x = y,\;\; \bar y = z,\;\; \bar z =  B x + G(y,z),
\label{3dHMgen}
\end{equation}
which are called {\em three-dimensional generalized H\'enon maps}. Such maps have the constant Jacobian $B$ and are, in a sense, the simplest three-dimensional nonlinear maps. Because of this, their use for the purpose of searching and studying strange attractors, including homoclinic ones, is very convenient and beneficial in many respects. In this paper we will use them as the main examples illustrating our results.

In Section~\ref{sec:geom} we give a qualitative description of the discrete Lorenz attractors dealing more accent to
their geometric properties and a comparison with the classical Lorenz attractors of three-dimensional flows.

One of
remarkable properties of discrete homoclinic attractors, including the mentioned above discrete Lorenz attractors, consists in the fact that they can arise in various kinds of models as a result of
rather simple and universal bifurcation scenarios \cite{GGS12,GGKT14}.
Moreover, these scenarios can be freely observed in one parameter families starting with those parameter values when the attractor is simple, e.g. a stable fixed point.
In Section~\ref{sec:scen} we give a phenomenological description of such scenarios for the case of discrete Lorenz-like attractors of orientable three-dimensional maps (the nonorientable case is discussed in Section~\ref{sec:othertype}). We show two typical examples of realizations of these phenomenological scenarios in one parameter families of concrete three-dimensional maps:  the three-dimensional H\'enon map (\ref{3dHM1}) and a Poincar\'e map for a Celtic stone model (see Figs.~\ref{3DHM1_evattr} and \ref{Kelt_lor}, resp.).
For comparison, we also show schematically in Fig.~\ref{Fig-scen1} the well-known scenario, see e.g. \cite{Sh80}, of the appearance of the Lorenz attractor in the one parameter family of the Lorenz model
\begin{equation}
\dot x = - \sigma (x-y), \; \dot y = - xz + r x - y, \; \dot z = xy - b z,
\label{eq:Lormod}
\end{equation}
where the parameters $b$ and $\sigma$ are fixed, $b=8/3, \sigma =10$, and the parameter $r$ is governing.

The fact that discrete Lorenz-like attractors can appear as a result of very simple and realistic bifurcation scenarios shows that they
should be found  in various models from applications.
As we know, the first such applied model was a nonholonomic model of Celtic stone \cite{GGK13}. More recently, discrete Lorenz attractors were also found in models of diffusion convection \cite{EM18,EM18a}.

In Section~\ref{sec:pseudo} we discuss one more, perhaps the most interesting fundamental property of discrete Lorenz attractors related to the fact that they can be genuine strange attractors,
i.e. robustly preserving their chaoticity at perturbations.
Recall that, as is commonly believed, a strange attractor is considered genuine,
if {\em all} its orbits have the positive maximal Lyapunov exponent and this property is open, i.e. holds for all close systems (in $C^r$-topology with $r\geq 1$).

Till recently, only hyperbolic attractors and flow Lorenz-like attractors (the latter belong to the class of singular hyperbolic attractors \cite{MPP98,S09})
%
%
%
could reliably be considered as genuine strange attractors.
However, in the paper \cite{TS98} by Turaev and Shilnikov, one more class of genuine strange attractors, the so-called {\em pseudohyperbolic attractors}, was introduced and, moreover,  an  example of such attractor for  a four-dimensional flow  was constructed.
%
This attractor was called in \cite{TS98} wild spiral attractor, since it contains a saddle-focus equilibrium together with
wild hyperbolic subsets and, hence, it allows homoclinic tangencies.\footnote{The wild hyperbolic sets were discovered by S.Newhouse \cite{N70,N79}. He used the term ``wild'' for notation of uniformly hyperbolic sets whose stable and unstable invariant sets (manifolds) have nontransversal intersections and this property holds for all close systems. Thus, such systems compose open regions (the so-called Newhouse regions)  in the space of dynamical systems and systems with homoclinic tangencies are dense in these regions.}
Despite the fact that such attractors were predicted (as well as geometrically constructed) in the late 90s, the first example of such an attractor in a
system of four differential equations was found much recently in \cite{GKT18}, see also \cite{GKT20}.

Discrete Lorenz attractors, e.g. those discovered in \cite{GOST05}, consist one more  class of wild pseudohyperbolic attractors in the case of three-dimensional maps.
In Section~\ref{sec:pseudo} we give a definition of pseudohyperbolic invariant sets of multidimensional diffeomorphisms (Def.~1) and explain how it
can be adapted to the discrete Lorenz attractors of three-dimensional maps. First, we discuss a local aspects of the corresponding theory related to the fact that such attractors can be born under codimension-3 bifurcations of fixed points. This allows to say about ``existence of genuine discrete Lorenz attractors''. Second, we review numerical methods for verification of pseudohyperbolicity of attractors when the local theory is not applicable. In particular, we demonstrate results of such verification for attractors from Figs.~\ref{fig:attrmap1}a,b.

In Section~\ref{sec:othertype} we consider other types of discrete Lorenz-like attractors including those observed in nonorientable three-dimensional maps and in Poincar\'e maps of periodically perturbed flows with the Lorenz attractors. We also discuss some questions related to both crises of discrete Lorenz-like attractors and creation  of the so-called discrete Lorenz-shape attractors of various types.
Main results of Section~\ref{sec:othertype} are new. In particular, we provide one more type of bifurcation scenario leading to the appearance of the so-called period-2 Lorenz attractor and give an example of its realization in the case of map (\ref{3dHM1}) with $B=-0.8$, see Fig.~\ref{Lor_per2_bif}. We study dynamical and geometrical properties of period-2 Lorenz attractors and show that these attractors can undergo crises leading to the emergence of strange homoclinic and heteroclinic attractors of new types. We construct geometric (skeleton) schemes of these attractors and their phase portraits in map (\ref{3dHM1}) with $B=-0.8$, see Figs.~\ref{Lor2_skel_p_phyp},~\ref{attr_skp_hom} and~\ref{attr_skp_het}. We show also that this period-2 Lorenz-like attractor
is pseudohyperbolic. On the other hand, the found us homoclinic attractors, containing the fixed point $O$ and a period-2 cycle, are not pseudohyperbolic because of the point $O$ is a saddle-focus. Moreover this point has a pair of complex conjugate multipliers close to 1:4 resonance (i.e. $\lambda_{1,2} = \rho e^{\pm i\varphi}$, where $\varphi$ is close to $\pi/2$ and $\rho<1$ is close to 1).
Thus, we also, indirectly, touched the problem of the structure of emerging homoclinic attractors when passing near strong resonances, which was formulated in \cite{GGS12}.

In Conclusion (Section~\ref{sec:concl}), we discuss a serious of open and future problems associated with discrete homoclinic attractors, emphasizing especial  attention to those related to the study of discrete Lorenz-like attractors.

\section{On geometrical properties of discrete Lorenz attractors of three-dimensional diffeomorphisms} \label{sec:geom}

By definition, discrete homoclinic attractors are strange attractors containing a single saddle fixed point of a map. As such an attractor is the stable closed invariant set, then it must contain also the closure of the unstable invariant manifold of the fixed point.

One of the main characteristics of discrete homoclinic attractors is the topological type of its fixed point $O$. First of all, we define the type $(i,j) = (\dim W^s(O), \dim W^u(O))$ of the point. Thus, in the three-dimensional case, hyperbolic fixed points may be of four types: (3,0) -- stable (sinks), (0,3) -- completely unstable (sources), and saddle fixed points of types (2,1) and (1,2).
Besides, saddle fixed points are divided into {\em saddles}, when all their three multipliers are real, and {\em saddle-foci}, when a pair of complex conjugate multipliers exists. In the last case, the corresponding discrete homoclinic attractors are usually called spiral.
According to [2,5], we divide these attractors into two types: discrete figure-8 spiral attractors when they contain a saddle focus (2.1), and discrete Shilnikov attractors when they contain a saddle focus (1,2).\footnote{Such discrete spiral attractors are often met in various applications, see e.g. \cite{BKS16,GKS17,GGKKB19,GSKK19}. Concerning the discrete Shilnikov attractors, they are also interesting from that point of view that they can give simple criteria of hyperchaos \cite{GSKK19}.}
%
%
%
In the case when $O$ is a saddle (2,1), we select, first of all, discrete homoclinic attractors of two types: Lorenz-like and figure-8 ones \cite{GGS12}. Discrete homoclinic attractors with saddles (1,2) also exist, but they are yet very poorly studied, some their examples see in \cite{GG16}.

As we said before, discrete Lorenz attractors were first found in \cite{GOST05}, see the examples
in Fig.~\ref{fig:attrmap1}a,b. Attractors in these figures contain a saddle fixed point $O$ of type (2,1) with multipliers $\lambda_1,\lambda_2,\gamma$ satisfying the following conditions
\begin{equation}
\begin{array}{l}
(a)\; 0< \lambda_1 <1,\;\; -1  < \lambda_2 <0, \;\; \gamma <-1,\;\; |\lambda_1\lambda_2\gamma|<1, \\
(b)\;\; |\lambda_1| > |\lambda_2|, \\
(c)\;\; \sigma = |\lambda_1\gamma| >1
\end{array}
\label{mults}
\end{equation}
The quantity $\sigma$ is called the saddle value and
it is defined, for a hyperbolic saddle periodic orbit, as the absolute value of product of two nearest to the unit circle stable and unstable multipliers, i.e. such multipliers that less and greater than 1 in modulus, respectively.


\begin{figure}[ht]
\centerline{\epsfig{file=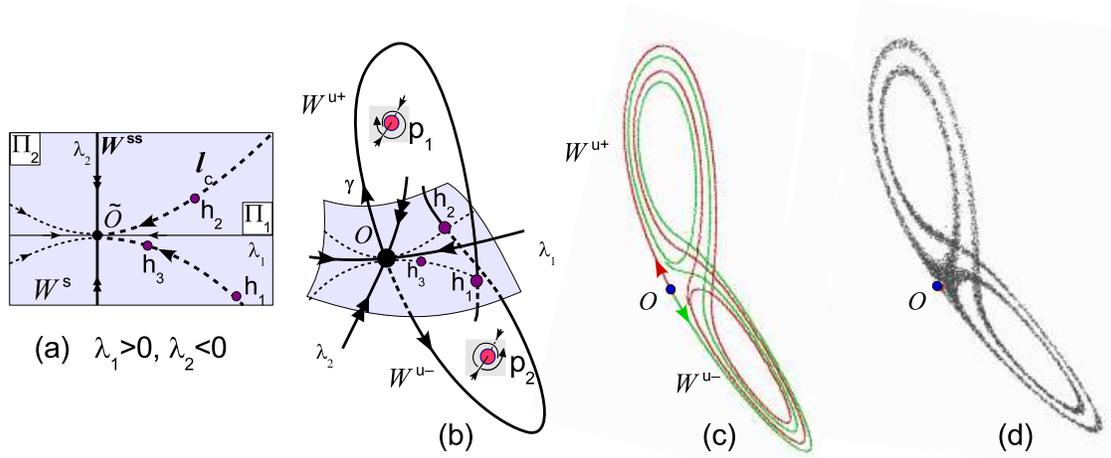, width=16cm
}}
\caption{{\footnotesize (a) An illustration for orbits behavior near
a  nonorientable sink $\tilde O$ having multipliers $-1<\lambda_2<0<\lambda_1<1$, where $|\lambda_2|<|\lambda_1|$; (b) a schematic homoclinic butterfly configuration of semi-global pieces of one-dimensional unstable invariant manifolds of the saddle point $O$ of type (2,1). In Figs. (c) and (d) ``live pictures'' for the case of map (\ref{3dHM1}) with $B=0.7; M_1=0; M_2=0.815$ are shown: (c)
the separatrices  $W^{u+}$ and $W^{u-}$ of the saddle $O$;
(d) a plot of the homoclinic attractor in an appropriate two-dimensional projection. }}
\label{DiscrLor}
\end{figure}

By virtue of (\ref{mults}), the unstable invariant manifold $W^u(O)$ of the point $O$ is one-dimensional and, hence, it is divided by the point $O$ into two connected components -- separatrices $W^{u+}$ and $W^{u-}$. Since the unstable multiplier $\gamma$ of $O$ is negative, $\gamma<-1$,  the separatrices $W^{u+}$ and $W^{u-}$ are invariant under $T^2$ and such that $T(W^{u+}) = W^{u-}$ and $T(W^{u-}) = W^{u+}$. This implies that points of $W^{u}$ jump under iterations of $T$ alternately from one separatrix to another.
The stable invariant manifold $W^s(O)$ of the point $O$ is two-dimensional, and the restriction of the map $T$ into the local stable manifold $W^s_{loc}(O)$,
i.e. the map
$T_s = T\bigl|_{W^s_{loc}}$, has a fixed point $\tilde O = O\cap W^s_{loc}$ that is a nonorientable node. By virtue of (\ref{mults}a) and (\ref{mults}b), in  $W^s_{loc}$ there is the so-called strong stable invariant manifold $W^{ss}$ that is a curve touching at $\tilde O$ to the eigendirection corresponding to the strong stable multiplier $\lambda_2<0$. Note that the curve $W^{ss}$ divides $W^s_{loc}$ into two parts $\Pi_1$ and $\Pi_2$, and, since
$\lambda_1>0$, points from $\Pi_1$ can not jump to $\Pi_2$ (and otherwise) under iterations of $T_s$, see Fig.~\ref{DiscrLor}a.

The two-dimensional map $T_s$ can be viewed, for simplicity, as a linear map $\bar x = \lambda_1 x, \bar y = \lambda_2 y$, where $0< \lambda_1 <1, -1  < \lambda_2 <0$ and $|\lambda_2|<|\lambda_1|$. Then $W^{ss}$ has the equation $x=0$, and
$\Pi_1 = \{(x,y) |x>0 \}$, $\Pi_2 = \{(x,y) |x<0 \}$. Moreover, as well-known, a neighborhood of the point $O$ is foliated into invariant (under $T_s^2$) curves $\{l_c\}$ of form $y = c |x|^\alpha$, where
$\alpha = \ln |\lambda_2|/\ln |\lambda_1|$. All these curves, except for the curve $x=0$, 
enter to $O$ touching the axis $y=0$. Let a point $h_1(x_1,y_1)$ belong to some curve $y = c x^\alpha$ with $x>0$ (or
$y = c (-x)^\alpha$ with $x<0$). Then its image, the point $h_2 = T_s(h_1) = (\lambda_1 x_1, \lambda_2 y_1)$, will belong to the curve $y = - c x^\alpha$ (resp.,
$y = - c (-x)^\alpha$). If $c\neq 0$, these two curves compose the boundary of
an exponentially narrow wedge adjacent to the point $O$ on one side (e.g. on $\Pi_1 $ as in Fig.~\ref{DiscrLor}a -- here the wedge is given by the inequalities
$|y|\leq c_1 x^\alpha, x\geq 0$, where the constant $c_1$ is defined by the initial data
$x_1$ and $y_1$, namely, $c_1 = y_1 x_1^{-\alpha}$).
Evidently, forward iterations under $T_s$ of the point $h_1$, i.e. the points $h_1,h_2,...$, where $h_{i+1}=T_s(h_i)$, will jump alternately on the sides of this wedge and will accumulate to the point $O$ as $i\to + \infty$.

Now one can imagine that $h_1$ is an intersection point of $W^{u+}$ with $W^s$. Then $h_2$ is the intersection point of $W^{u-}$ with $W^s$, since $T(W^{u+}) = W^{u-}$, and $h_3$ is again
%
%
an intersection point of $W^{u+}$ with $W^s$, etc. Correspondingly, the points $h_1,h_2,...$ are points of some homoclinic to $O$ orbit, i.e. an orbit entirely lying in the invariant set $W^u(O)\cap W^s(O)$ and distinct of $O$. We call the points $h_1,h_2,...$ {\em primary homoclinic points}. Evidently, their disposition in $W^s_{loc}$ defines in much a skeleton of a homoclinic attractor that, in turn, can be viewed as the closure of the unstable manifold $W^u(O)$.

It is easy to see that
the configuration of semi-global pieces of $W^{u+}$ and $W^{u-}$,
till the first their intersections with $W^s_{loc}$, see Fig.~\ref{DiscrLor}b, will be similar to the configuration of homoclinic butterfly in the Lorenz model (\ref{eq:Lormod}),
see Fig.~\ref{Fig-scen1}b.
Especially this similarity is noticeable when we consider longer pieces of separatrices, Fig.~\ref{DiscrLor}c, or take many iterations of a single point on the attractor, Fig.~\ref{DiscrLor}d.
Therefore,
it is no coincidence that it often turns out that the discrete Lorenz attractors found look very similar to the classical Lorenz attractor.



\section{On phenomenological scenarios of onset of discrete Lorenz attractors.} \label{sec:scen}

\begin{figure}[ht]
\centerline{\epsfig{file=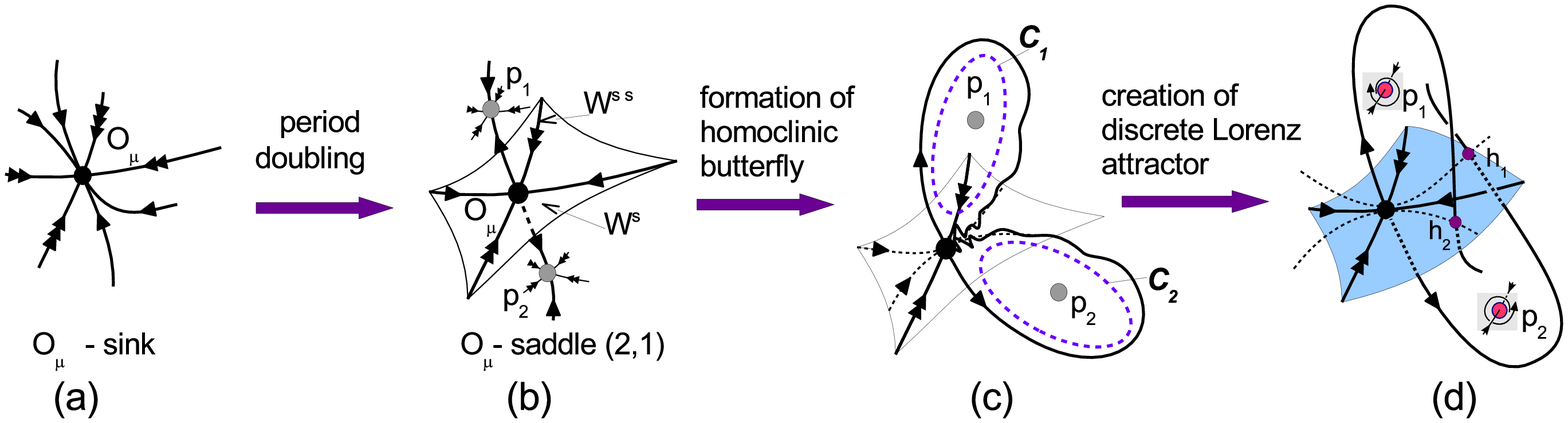, width=18cm
}}
\caption{\footnotesize{A phenomenological description of main steps for the scenario of onset of the discrete Lorenz attractor: (a) a stable fixed point $O_\mu$ $\rightarrow$ (b) after a period-doubling bifurcation: the point $O_\mu$ becomes a saddle (2,1) and the attractor now is  a period-2 cycle $(p_1,p_2)$  $\rightarrow$ (c) a creation of a homoclinic butterfly configuration of $W^u(O_\mu)$ and a formation of a closed invariant curve
$(C_1,C_2)$ of period 2, i.e. $T_\mu(C_1)=C_2$ and $T_\mu(C_2)=C_1$ (this curve is saddle if $\sigma>1$) $\rightarrow$ (d) as a variant (see option [sc1] below), the cycle $(p_1,p_2)$ loses stability merging with the saddle curve $(C_1,C_2)$ and the discrete Lorenz attractor appears as the unique attracting invariant set.  }}

\label{scen-Lor1}
\end{figure}

One of the most interesting peculiarities of discrete homoclinic attractors is that they can occur in one parameter families of (three-dimensional) maps
as a result of fairly simple universal bifurcation scenarios
\cite{GGS12,GGKT14,GG16}.\footnote{The idea of constructing such scenarios goes back to the paper \cite{Sh86} by Shilnikov, in which he described a phenological scenario of spiral chaos appearance in the case of multidimensional flows.} This fully applies to the discrete Lorenz attractors. Schematically, main stages of the corresponding scenarios
are shown in Fig.~\ref{scen-Lor1} for the case of a one parameter family $T_\mu$ of three-dimensional orientable maps.

This scenario begins from those values of $\mu$ where the map $T_\mu$ has
an asymptotically stable fixed point (sink) $O_\mu$. We assume that at  a certain value of $\mu$
the point $O_\mu$ loses its stability under the {\em supercritical period-doubling bifurcation}: $O_\mu$ becomes a saddle (2,1) and a stable cycle
$(p_1, p_2)$ of period 2 is born
(at this moment the cycle becomes the attractor).
When the parameter changes,  this cycle and all attracting invariant sets generated from it can lose stability as a result of a set of certain bifurcations.
How this happens depends on specific of the problem.
We select two the simplest and most common options:

\begin{itemize}

\item[{\bf[sc1]}]
the cycle $(p_1, p_2)$ loses the stability under the {\em subcritical} Andronov-Hopf bifurcation for maps: a period-2 closed  invariant curve $(C_1, C_2)$ of saddle type merges with the stable cycle $(p_1, p_2)$ and the latter becomes saddle (of type (1,2));

\item[{\bf[sc2]}]
the stable cycle
$(p_1,p_2)$ undergoes the {\em supercritical} Andronov-Hopf bifurcation for maps, after which the cycle becomes saddle (of type (1,2)) and a period-2 stable closed invariant curve $(S_1,S_2)$ is born, and then the stable curve $(S_1,S_2)$ merges with the saddle curve
$(C_1, C_2)$ and both
disappear.
\end{itemize}

Note that the above period-2 invariant curve
$(C_1, C_2)$ is formed from the homoclinic butterfly configuration of unstable separatrices $W^{u+}$ and $W^{u-}$ of $O_\mu$ at that moment
\begin{figure}[ht]
\centerline{\epsfig{file=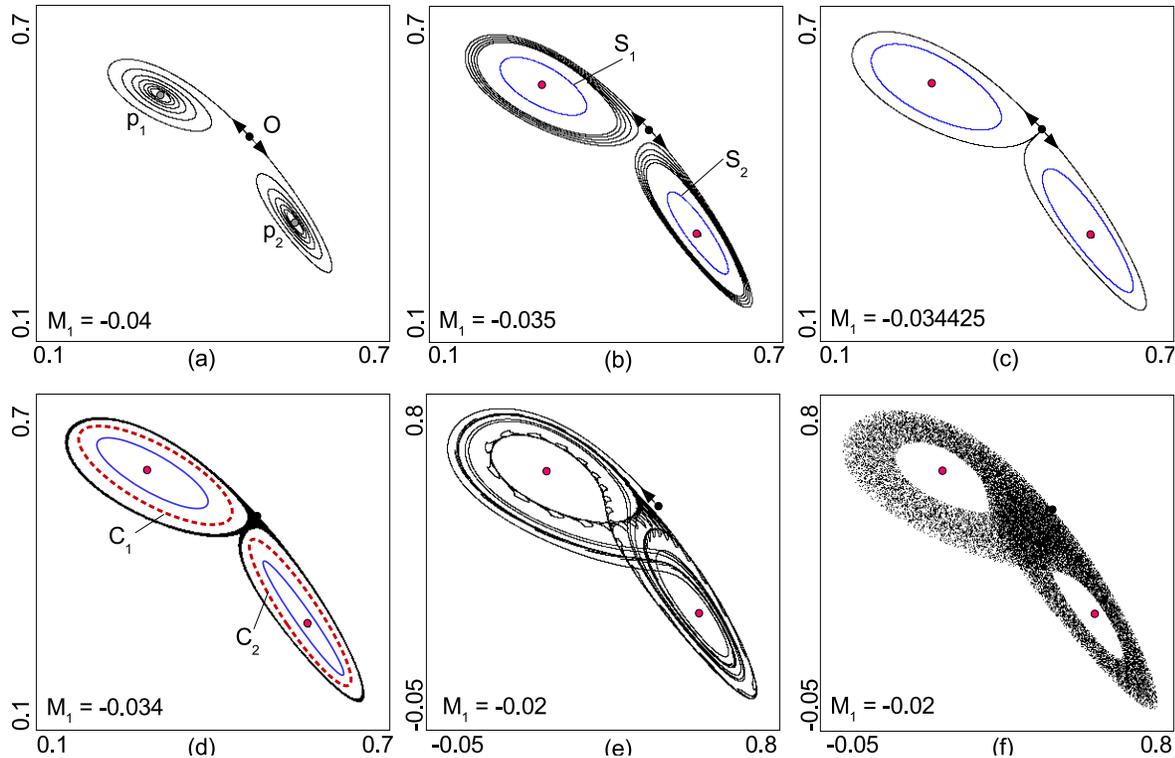, width=16cm
}}
\caption{{\footnotesize
Evolution of attractors in map (\ref{3dHM1}) with $B = 0.7, M_2 = 0.85$ when varying $M_1$, here the option {\bf[sc2]} of the scenario is realized:  (a) attractor is a period-2 cycle
$(p_1,p_2)$; (b) attractor is a stable closed invariant curve $(S_1,S_2)$ of period 2 that was born from the cycle
$(p_1,p_2)$); (c) creation of a homoclinic butterfly of the saddle $O$; (d) a saddle closed invariant curve $(C_1,C_2)$ is formed from the butterfly; (e)--(f) the curves  $(S_1,S_2)$ and $(C_1,C_2)$ merge and disappear and the discrete Lorenz-like attractor becomes unique attracting set: in (e) a behavior of one unstable separatrice of $O$ is shown (other separatrice behaves symmetrically), in (f) the phase portrait of attractor is shown (here $O=(x,x,x)$, where $x\approx 0.51$).  }}
\label{3DHM1_evattr}
\end{figure}
when the invariant manifolds of $O_\mu$ begin to intersect, see Fig.~\ref{scen-Lor1}c and Fig.~\ref{3DHM1_evattr}d. Besides, the period-2 curve $(C_1, C_2)$ is saddle when $\sigma>1$ and stable when $\sigma <1$.
\begin{figure}[!ht]
\centerline{\epsfig{file=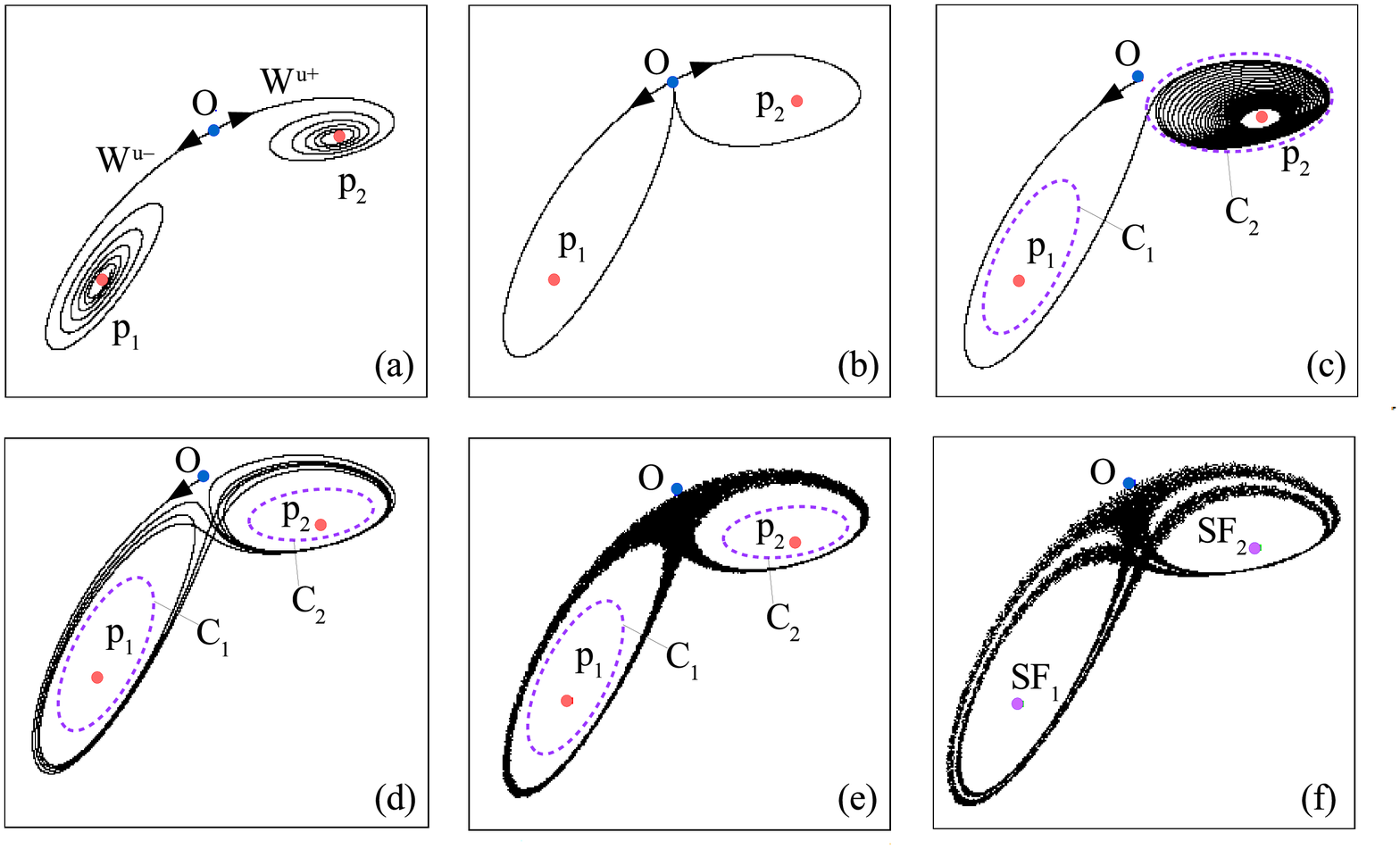, width=16cm
}}
\caption{{\footnotesize
Evolution of attractors in a nonholonomic model of Celtic stone, based on the paper \cite{AGJS19},  when varying a parameter $E$ (the full energy of the stone):(a) $E=741$,  attractor is a period-2 cycle $(p_1,p_2)$ (appearing after period doubling bifurcation with a stable fixed point); (b) $E=742.063$, creation of a homoclinic butterfly of the separatrices $W^{u+}$ and $W^{u-}$ together with the saddle $O$; (c) $E=742.2$, the separatrices are reconstructed and a saddle closed curve $(C_1,C_2)$ of period 2 is formed; (d)--(e) $E=742.5$, multistability -- the stable cycle  $(p_1,p_2)$ and a discrete Lorenz-like attractor coexist, in (d) and (e) there are shown the separatrix $W^{u+}$ ($W^{u-}$ is symmetric) and  the phase portrait, respectively; (f) $E=743.3$, the curves $C_1$ and $C_2$ merge with the cycle  $(p_1,p_2)$ and the discrete Lorenz attractor becomes the unique attracting invariant set.  Here the option {\bf[sc1]} of the scenario is realized. }}
\label{Kelt_lor}
\end{figure}


It is worth to mention that both variants, {\bf[sc1]} and {\bf[sc2]},  are often met in applications. In particular, we found their examples in three-dimensional generalized H\'enon maps of form (\ref{3dHMgen}), see \cite{GOST05,GGS12,GGOT13},
when the value $0<B<1$ of the Jacobian is not too small\footnote{When $B$ is close enough to zero (experimentally, when $B$ is less than $0.3$ or $0.4$ and so), even if the conditions (\ref{mults}) hold, the cycle $(p_1,p_2)$ has a tendency to lose its stability via a cascade (finite or infinite) of period-doubling bifurcations. Then a form of resulting attractor can look as ``severely deformed Lorenz-like'' one or as ``thickened H\'enon-like'' one. That is, a certain analogy between phase portraits of discrete and flow Lorenz attractors can be not clearly observed and even fully disappear.},  and in a model of Celtic stone \cite{GGK13,AGJS19}. Two illustrations of the scenarios are shown in Fig.~\ref{3DHM1_evattr} (for the case of three-dimensional
H\'enon map (\ref{3dHM1}))  and in Fig.~\ref{Kelt_lor} (for the case of a nonholonomic model of Celtic stone, \cite{GGK13,AGJS19}) where variants {\bf[sc2]} and {\bf[sc1]} of the general scenario are observed, respectively.

\begin{figure}[ht]
\centerline{\epsfig{file=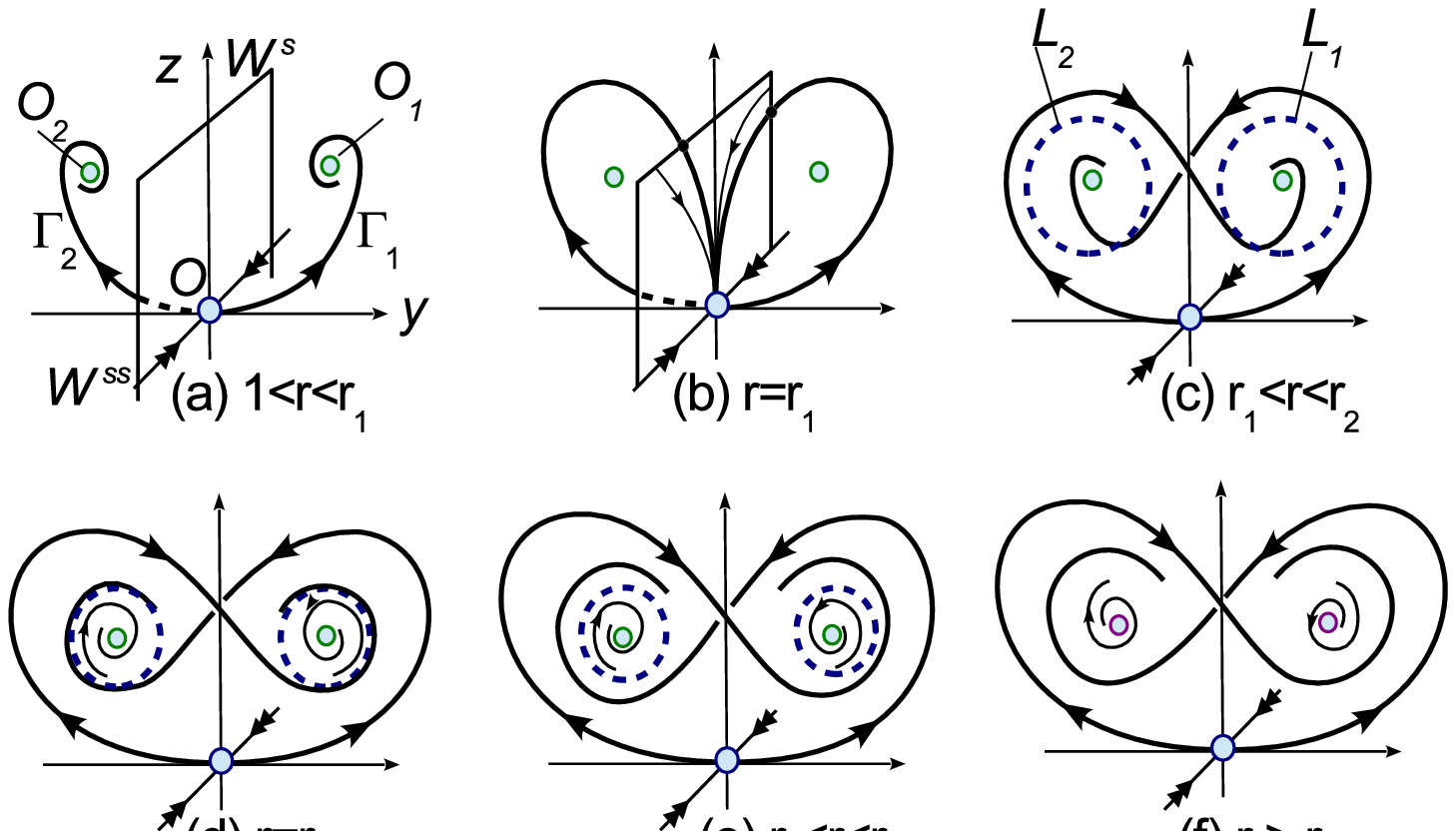, width=16cm
}}
\caption{{\footnotesize Scenario of onset of the Lorenz attractor in the Lorenz model (\ref{eq:Lormod}) with $b=8/3, \sigma = 10$ and varying $r$,  based on the paper \cite{Sh80}. Here $r_1 \approx 13.92$,  $r_2 \approx 24.06$ and $r_3 \approx 25.06$. (a) the attractors are non-zero equilibria
$O_1$ and $O_2$; (b) the homoclinic-butterfly is created; (c) saddle limit cycles $L_1$ and
$L_2$ are born; (d) the separatrices $\Gamma_1$ and $\Gamma_2$ belong to the stable manifolds of the limit cycles $L_2$ and
$L_1$, respectively; (e) multistability: the Lorenz attractor coexists with the attractors $O_1$ and $O_2$; (f) $O_1$ and $O_2$ simultaneously lose stability under a  subcritical Andronov-Hopf bifurcation for flows
and the Lorenz attractor becomes unique attractor of the model. }}
\label{Fig-scen1}
\end{figure}

Note that the scenarios with options {\bf[sc1]} or {\bf[sc2]} are quite expected here, because they are, in a sense, inherited from flows.
So, in Fig.~\ref{Fig-scen1}, a schematic picture is shown for the scenario of onset of the Lorenz attractor in the Lorenz model (\ref{eq:Lormod}) when varying the parameter $r$. Here two saddle limit cycles symmetrical to each other are born
from a homoclinic butterfly and then these cycles merge with the non-zero stable equilibria, and, just after this subcritical Andronov-Hopf bifurcation, the Lorenz attractor appears to be the unique attracting invariant set of the model. Thus, the flow analogue of option {\bf[sc1]} is realized here.

The same scenario can be also observed in the Shimizu-Morioka model
\begin{equation}
\dot X = Y, \; \dot Y = X(1-Z) -\lambda Y, \; \dot Z = -\alpha Z - X^2.
\label{eq:SMsys}
\end{equation}
However, the flow analogue of option {\bf[sc2]} can also be observed here \cite{ASh86,ASh93}:  the non-zero equilibria lose their stability softly and  simultaneously under a supercritical Andronov-Hopf  bifurcation and, as result, two stable symmetric each other limit cycles are born and, further, these cycles merge with saddle ones and both disappear.

In general, the stage of the scenario when global bifurcation occurs associated with the formation of homoclinic structures at the point $O_\mu $ (i.e., when the invariant manifolds
$ W^u(O_\mu) $ and $W^s(O_\mu)$ begin to intersect) is key one.
It indicates the appearance of comp\-li\-ca\-ted dynamics of the map, and, if
the unstable invariant manifold
$W^u(O_\mu)$ entirely lies in the absorbing region, then  a discrete homoclinic attractor can emerge.
Initially, this attractor can coexist with other stable regimes (e.g. with the stable cycle
$(p_1,p_2)$, see Fig.~\ref{Kelt_lor}c.
However, when the latter have lost the stability, the homoclinic attractor remains the unique attracting set of the map.

As we explained before, when conditions (\ref{mults}) hold,  the corresponding discrete homoclinic attractor will
look like the Lorenz attractor. However,
we are dealing here with its discrete version, when the fixed point $O_\mu$ plays a role of the zero equilibrium of the Lorenz system (\ref{eq:Lormod}) and
the saddle cycle $(p_1, p_2)$ of period 2 resides inside the ``holes'', instead of two non-zero equilibria $O_1$ and $O_2$.
In addition,
the one-dimensional unstable manifold
$W^u(O_\mu)$ consists of two connected components, separatrices $W^{u+}$ and $W^{u-}$, and, since the unstable multiplier
$\gamma$ of $O_\mu$ is negative, $\gamma < -1$, the points on $W^u(O_\mu)$ will ``jump''  from one separatrix to other, whereas,
in the classical Lorenz attractor, each of the separatrices of zero equilibrium is invariant itself. Note also that the negativity of multipliers of the fixed point in the discrete attractor provides a local (on its invariant manifolds) symmetry which, in fact, plays the role of
the global symmetry $x\to -x, y\to -y, z\to z$ in the Lorenz model (\ref{eq:Lormod}).

\section{Pseudohyperbolicity of discrete Lorenz attractors.} \label{sec:pseudo}

It is worth to mention that the discrete Lorenz attractor, the same as its flow analogue, can be genuine. We note at once that the condition $\sigma >1$ is simply necessary for this event. First, when $\sigma <1$ a homoclinic attractor may not arise at all. Here, the homoclinic-butterfly bifurcation leads to the birth of period-2 invariant curve $(C_1,C_2)$ that will be stable, i.e. this curve becomes attractor. Moreover, its further evolution can lead to attractors of torus-chaos type (see e.g. \cite{AfrSh83}) or so, but not Lorenz-like ones.

However, there is another much more delicate and interesting reason. Even if a discrete Lorenz-like attractor arises, it will always be the {\em quasiattractor} in the case $\sigma<1$.
That is, it will be an attractor that looks strange but only from ``physical point of view'' \cite{AfrSh83b}. The fact
is that quasiattractors, unlike genuine strange attractors, either contain stable periodic orbits (sinks)
with very narrow domains of attraction, or such sinks appear at arbitrary small perturbations.
Usually, such stable orbits are not detected in experiments and the attractor seems chaotic, but, sometimes, these orbits come to light for certain values of parameters called windows of stability.

One of the main reasons for such nonrobust chaoticity of homoclinic quasiattractors
is the inevitable appearance of homoclinic tangencies, i.e. nontransversal intersections of stable and unstable invariant manifolds of the saddle fixed point.
When $\sigma <1$, bifurcations of such tangencies lead to the birth of stable periodic orbits \cite{GaS73,G83,GST97a}.
The latter have,
as rule, very narrow domains of attraction and big periods.
However, everything can change in the case $\sigma> 1$. Here, in general,  homoclinic tangencies do not destroy the chaoticity (their bifurcations lead to the birth of saddle periodic orbits \cite{GST93c,GST08}, instead stable ones at $\sigma<1$). Therefore
the discrete Lorenz-like attractors can keep their chaoticity when changing parameters.
Such robustness of chaos is a direct consequence of the fact that the discrete Lorenz-like attractors
can be pseudohyperbolic.

Speaking shortly,  an attractor is {\em pseudohyperbolic} if, in some its neighborhood $D$ (one can think that $D$ is some absorbing region), a ``weak'' version of hyperbolicity is fulfilled, i.e. there is an exponential contraction along some invariant directions and exponential expansion of volumes on the subspaces transversal to them, for more details see \cite{GKT18arx,GGKK18}.

In general setting,  the problem of searching for genuine (pseudohyperbolic) attractors seems to be wide open. Especially great difficulties arise here in that part of the problem where it is required to distinguish attractors found from quasiattractors.
Only for some types of such attractors, their pseudohyperbolicity has been proved. In particular, this relates to the discrete Lorenz attractors, see e.g. \cite{GOST05,GGKK18}.
Moreover, it seems that the existence of discrete Lorenz-like attractors in three-dimensional maps should be quite common phenomenon.
This, in particular, is due to the fact that such attractors can arise as a result of rather simple global scenarios described in Section~\ref{sec:scen}, as well as a result of local
bifurcations, which we discuss below.




\subsection{Discrete Lorenz-like attractors under local bifurcations.}

In the paper \cite{SST93} it was shown that the flow Lorenz attractors can be born as result of local bifurcations of an equilibrium with three zero eigenvalues. In the same paper it was noted that, in the case of three-dimensional maps, local bifurcations of a fixed point with a triplet $(-1,-1,+1)$ of multipliers can lead to the birth of discrete Lorenz-like attractor. In this case, the second power of the map may be locally (near the fixed point) represented as a map $o(\tau)$-close to the time-$\tau$ map of the Shimizu-Morioka system (\ref{eq:SMsys}), where $\tau$ can be taken as small as we want. Since system (\ref{eq:SMsys}) has the Lorenz attractor for some open domain of positive parameters $(\alpha,\lambda)$ \cite{ASh86,ASh93,OT17}, it follows that the original map has an attractor that, if to take every second iteration,  is a $\tau$-periodic perturbation of a Lorenz attractor.

Map (\ref{3dHM1}) turned out to be the first concrete model in which such codemension 3 bifurcation was studied \cite{GOST05}.
Note that map (\ref{3dHM1}) has the fixed point with a triplet $(-1,-1,+1)$ of multipliers for the values $A^* = (M_1 = - 1/4, M_2 = 1, B = 1)$ of the parameters. It was shown in \cite{GOST05} that a small discrete Lorenz attractor exists in (\ref{3dHM1}) for some domain of the parameters adjoining to $A^*$ from the half-space $B<1$ (see also Proposition 1 below).

In the paper \cite{GGOT13} this result was extended to the case of three-dimensional generalized H\'enon maps (\ref{3dHMgen}). When map (\ref{3dHMgen}) has a fixed point, this point can be moved to the origin, and then the map takes the form
\begin{equation}
\bar x = y,\; \bar y = z, \; \bar z = B x + C y + A z + a y^2 + b yz + c y^2 + ...,
\label{eq:3mabc}
\end{equation}
where the dots stand for cubic and higher order terms. The characteristic equation at the fixed point $O(0,0,0)$ is
$$
\lambda^3 - A\lambda^2 - C\lambda - B = 0.
$$
Thus, the fixed point has multipliers $(+1,-1,-1)$ at $(A=-1,C=1,B=1)$ and, hence, for nearby values of the parameters, map (\ref{eq:3mabc}) can be written as
\begin{equation}
\bar x = y,\; \bar y = z, \; \bar z = (1-\varepsilon_1) x + (1-\varepsilon_2) y - (1+\varepsilon_3) z + a y^2 + b yz + c y^2 + ...
\label{eq:3mabc2}
\end{equation}

The following result was established in \cite{GGOT13} (see Lemma 3.1 there).\\

\textbf{Proposition 1.}  {\em Let the following condition hold}
\begin{equation}
(c-a)(a-b+c) > 0.
\label{eq:3mabc3}
\end{equation}
{\em Then map {\rm (\ref{eq:3mabc2})} has a \textsf{a discrete pseudohyperbolic Lorenz-like attractor} for all $\varepsilon$ from an open, adjoining to $\varepsilon =0$, subregion ${\cal D}$ of
$\{\varepsilon_1>0, \varepsilon_1 + \varepsilon_3 >0, |\varepsilon_2 -\varepsilon_1 - \varepsilon_3|\leq L(\varepsilon_1^2 + \varepsilon_3^2)\}$, for some $L>0$.} \\

We note that the inequality (\ref{eq:3mabc3}) satisfies automatically for map (\ref{3dHM1}), for which $c=-1,a=0,b=0$. However, for example, for the map
\begin{equation}
\bar x = y,\; \bar y = z, \; \bar z =M_1 + B x + M_2 z - y^2,
\label{eq:3mabc4}
\end{equation}
that is well-known as a ``homoclinic map'' \cite{GST93c,GGS10}, where $a=-1,c=0,b=0$, the inequality (\ref{eq:3mabc3}) is not fulfilled. In principle, this is not surprisedly, since the map (\ref{eq:3mabc4}) is inverse to the map (\ref{3dHM1}), and, hence, it should have a repeller instead the attractor.

Note that map (\ref{eq:3mabc4}) was introduced in the paper [26], where its main bifurcations were studied. When $B=0$ this map becomes effectively two-dimensional map of form $\bar y = z, \; \bar z = M_1 + M_2 z - y^2$. It is the well-known two-dimensional endomorphism introduced by C. Mir\'a \cite{Mira64} yet in 60s. Therefore we will call map (\ref{eq:3mabc4}) the three-dimensional Mir\'a map. As it was shown in \cite{GGS12,GGKT14}, map (9) demonstrates chaotic dynamics of a completely different type than the map (\ref{3dHM1}) -- it is more associated  with spiral attractors than with Lorenz-like ones.

\subsection{On numerical verification of pseudohyperbolicity.} \label{LMPmethod}

Let us return to the attractors of Figs.~\ref{fig:attrmap1}a,b, for which we now consider the question of their pseudohyperbolicity. Note that the corresponding values of the parameters ($ M_1 = 0, M_2 = 0.85, B = 0.7$ and $ M_1 = 0, M_2 = 0.815, B = 0.7 $) are not nearly close to $A^* = (M_1 = - 1/4, M_2 = 1, B = 1) $. Thus, we can not deduce the desired pseudohyperbolicity of these attractors from the local analysis, e.g. due to Proposition 1.
It should be verified additionally that is not quite simple exercise. However, one can make this numerically. First we note that the fixed point $O(x=y=z=1-M_2-B)$ of the attractor is pseudohyperbolic itself: it has multipliers satisfying condition (\ref{mults}). Thus, if the attractors are pseudohyperbolic, then $\dim N_1 =1, \dim N_2 =2$.   This follows also that the spectrum of Lyapunov exponents $\Lambda_1 > \Lambda_2 > \Lambda_3$ on the attractors must be such that the inequalities
\begin{equation}
\Lambda_1 > 0, \; \Lambda_1 +   \Lambda_2 > 0, \; \Lambda_1 + \Lambda_2 + \Lambda_3  <0
\label{lyap}
\end{equation}
are fulfilled.

\begin{figure}[ht]
\centerline{\epsfig{file=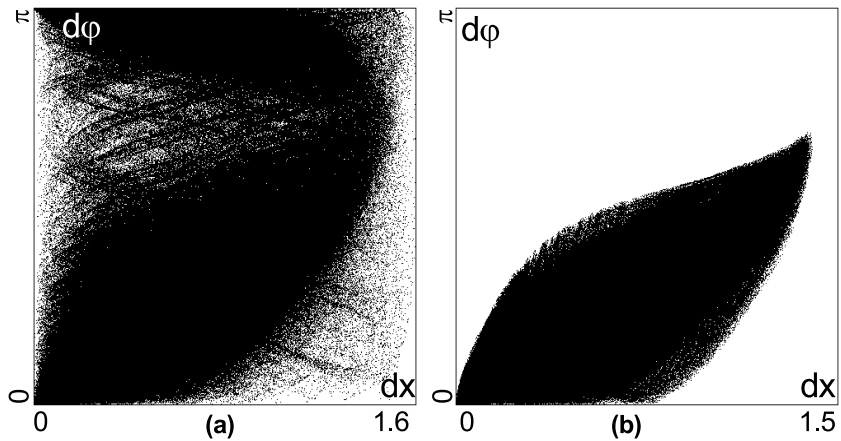, width=12cm}}
\caption{{\footnotesize LMP-graphs for attractors of (a) Fig.~\ref{fig:attrmap1}a  and (b) Fig.~\ref{fig:attrmap1}b.}}
\label{LMPdL2}
\end{figure}

It is shown standardly that these inequalities
are valid for numerically found exponents $\Lambda_i$ (for both the attractors of Figs.~\ref{fig:attrmap1}a,b).
However, such Lyapunov exponents are certain average characteristics for the orbits of attractor, therefore, in principle, it is possible that the attractor has very small ``gaps'' (whose sizes may be less than any reasonable accuracy of calculations), where conditions (\ref{lyap}) on exponents for the corresponding orbits from ``gaps'' are violated. In order to exclude
(as confident as possible)
this situation, one can additionally check numerically conditions from Definition~\ref{df:pseudo}. As it seems, the condition (ii) from this definition looks to be the most delicate and indicative. It says, in particular, that the field $N_1(x)$ of strongly contracting directions depends continuously on points $x$ of attractor. This field can be calculated in various ways. In particular, one of such methods, the so-called LMP-method (abbreviation for ``Light Method of Pseudohyperbolicity'' checking), was proposed in \cite{GKT18arx}, see also \cite{GGKK18}. This method allows to construct a field of vectors corresponding to the strongest contractions and to display, in form of LMP-graph, the dependence of the angles $d\varphi$ between vectors $N_1(x)$ and $N_1(x_2)$ on the distance $dx$ between the corresponding points $x_1$ and $x_2$ on the attractor.

In Fig.~\ref{LMPdL2} we represent the LMP-graphs for the discrete Lorenz attractors of Fig.~\ref{fig:attrmap1}a (left) and Fig.~\ref{fig:attrmap1}b (right).
The graph from Fig.~\ref{LMPdL2}b shows that the corresponding attractor is pseudohyperbolic with a lot of confidence. Indeed, here $d\varphi \to 0$ as $dx \to 0$, as it should be when the field $N_1$ is continuous. On the other hand, the graph from Fig.~\ref{LMPdL2}a
looks quite ``chaotic'' nearly the axis $dx=0$: it seems that if $dx\to 0$ then various sequences $d\varphi = d\varphi(dx)$ can take any partial limits in $[0,\pi]$ (in principle, as a tendency of calculations demonstrates, the LMP-graph approaches to a ``black square'' in this case). This shows that here the condition (ii) of continuity of field $N_1$ from Def.~\ref{df:pseudo} is not fulfilled.  The latter is characteristic for quasiattractors.

\section{Other types of discrete Lorenz-like attractors.} \label{sec:othertype}

One of the most natural ways to get discrete Lorenz-like attractors for three-dimensional maps was proposed in the paper \cite{TS08}, in which it was shown that such attractors can exist in Poincar\'e maps for periodically perturbed systems having the Lorenz attractor. When these perturbation are sufficiently small, the attractor is pseudohyperbolic.  Such discrete Lorenz attractor contains the unique fixed point, a saddle of type (2,1), with positive multipliers ($0<\lambda_1<\lambda_2<1<\gamma$), and, besides,  two fixed points of the map (both saddle-foci of type (2,1)) reside in two holes of the attractor. These characteristics of such attractors make them significantly different from the discrete Lorenz attractors considered in sections 2 and 3. Recall that the latter have fixed points with negative unstable multiplier and period-2 points reside in two its holes.

However, for the second power of map, these attractors may not be distinguishable at all.
Nevertheless, these attractors are different even in this case. In particular,
the discrete Lorenz attractor always possesses a symmetry that is inherited by the local symmetry between  $W^{u+}$ and $W^{u-}$ due to the unstable multiplier of the fixed point $O$ is negative.
When the unstable multiplier is positive, the behavior of separatrices $W^{u+}$ and $W^{u-}$ is independent and, thus, the corresponding discrete attractor will be not symmetric, in general.
This implies also that scenarios of creation of such attractors must differ from those considered in Section~\ref{sec:scen}, see Fig.~\ref{scen-Lor1}.\footnote{Here, instead of a codimension 1 period doubling bifurcation, one should have either a codimension 2 pitch-fork bifurcation (symmetric case), when the point $O$ becomes a saddle (2,1) and the fixed points $p_1$ and $p_2$ are born, see Fig.~\ref{scen-Lor1}b for the case where the points $O_\mu,p_1$ and $p_2$ are fixed , or such a configuration (two stable and one saddle fixed points) is created sharply due to a saddle node bifurcation, when the point $O$ remains stable and two new fixed points, saddle and stable, appear (in general, asymmetric, case).}

\begin{figure}[ht]
\centerline{\epsfig{file=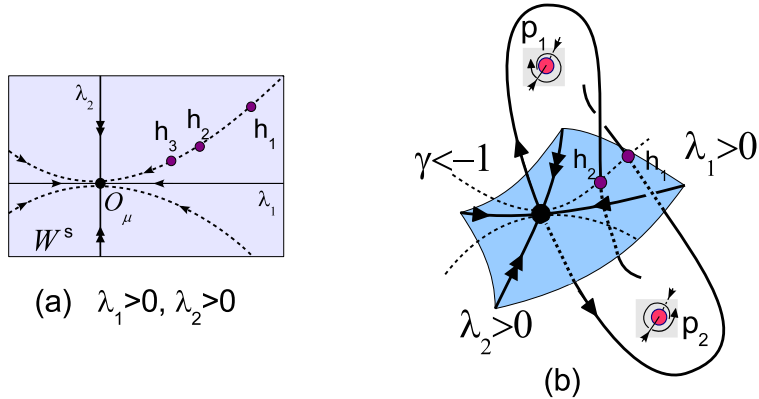, width=12cm
}}
\vspace{-1cm}
\caption{{\footnotesize (a) Behavior of iterations of points near a stable node with $0<\lambda_2<\lambda_1 <1$; (b) a homoclinic-butterfly configuration of semi-global pieces of unstable separatrices of a nonorientable saddle with multipliers $0<\lambda_2<\lambda_1<1$ and $\gamma < -1$.  }}
\label{Discr-Lor_nor}
\end{figure}

Discrete Lorenz-like attractors can also exist in three-dimensional {\em nonorientable} maps and, moreover, they can appear in one parameter families as result of scenarios similar to those in the orientable case.
Here only a certain difference in the structure of primary homoclinic orbits exists. It is shown in Figs.~\ref{DiscrLor}a
and~\ref{Discr-Lor_nor}a that homoclinic points $h_1,h_2,...$ in $W^s$  are posed differently: in the case of orientable attractor, Fig.~\ref{DiscrLor}a, they lie alternately on boundaries of
an exponentially narrow wedge, and in the case of nonorientable attractor, schematically presented in Fig.~\ref{Discr-Lor_nor}b, they lie on the same invariant curve (of the form $y = c |x|^\alpha$) entering to $O_\mu$, since the point $O_\mu$ has here the multipliers $0<\lambda_2<\lambda_1<1, \gamma < -1$.

\subsection{On period-2 Lorenz-like attractors}

As we discussed in Section~\ref{sec:scen} for the case of orientable maps, a good scenario leading to the emergence of the discrete Lorenz attractor must
contain bifurcation steps [s1] or [s2].  The same should be true for the nonorientable case.
However, good examples of implementations of such scenarios in the case of non-orientable maps have not yet been found (we think that their finding is only a matter of time).
Instead this, we found scenarios including several successive period doubling bifurcations with the cycle $(p_1,p_2)$. The simplest option here is as follows:

\begin{figure}[ht]
\centerline{\epsfig{file=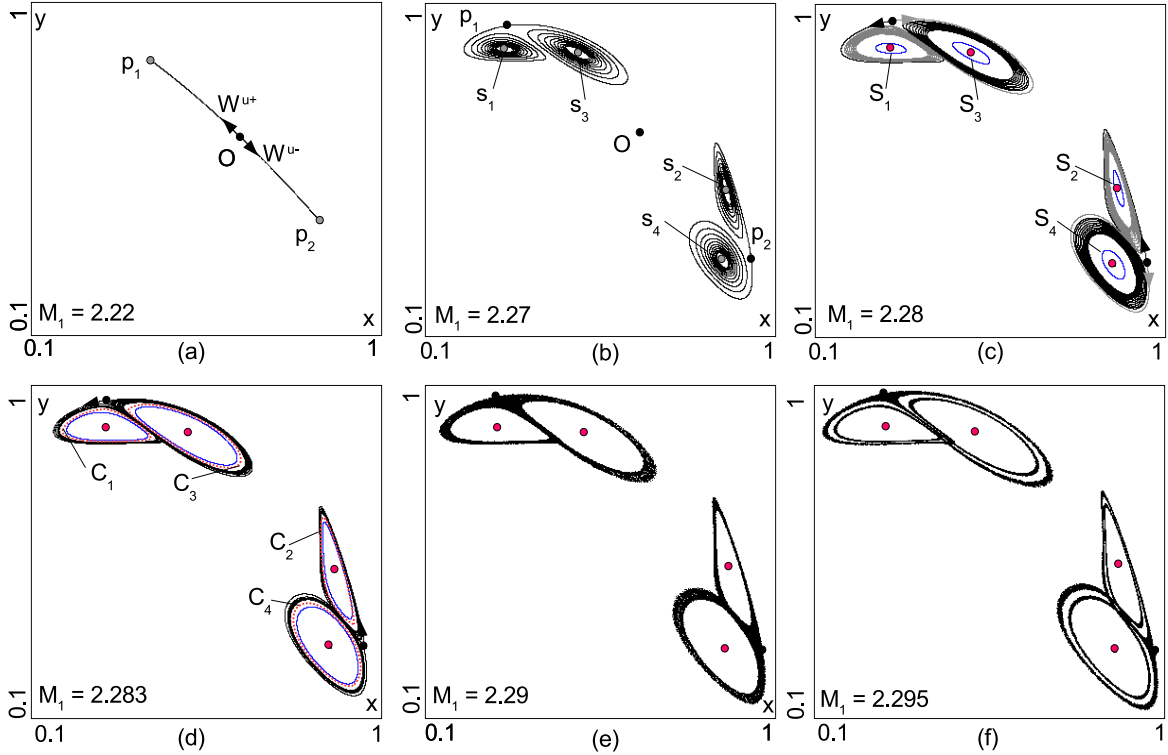, width=16cm
}}
\vspace{-1cm}
\caption{{\footnotesize Towards a scenario of onset of a discrete period-2 Lorenz-like attractor. The bottom row corresponds to strange attractors that contains the points $p_1(a,b,a)$ and $p_2(b,a,b)$,  where
(d) $a\approx 0.29, b\approx 0.96$; (e) $a\approx 0.28, b\approx 0.97$ and (f) $a\approx 0.27, b\approx 0.98$.
}}
\label{Lor_per2_bif}
\end{figure}

\begin{itemize}

\item[{\bf[sc3]}]
the cycle $(p_1,p_2)$ loses the stability under the second period doubling  bifurcation:
the cycle $(p_1,p_2)$ becomes saddle and a period-4 stable cycle $\tilde S =(s_1,s_2,s_3,s_4)$ is born, then the cycle $\tilde S$ loses stability (under an Andronov-Hopf bifurcation) according to
the options {\bf[sc1]} or {\bf[sc2]} for $T^2$.

\end{itemize}

Note that, {\em for the map $T^2$}, the points $p_1$ and $p_2$ become fixed and, thus,  the option {\bf[sc3]} gives a possibility for the simultaneous appearance of two discrete Lorenz-like attractors associated with these saddle fixed points.
This means that, {\em for the map $T$ itself}, a {\em period-2 Lorenz-like attractor} with four holes around points $s_1,s_2,s_3,s_4$  can arise.

In Fig.~\ref{Lor_per2_bif} we illustrate the {\bf[sc3]}-scenario for the nonorientable three-dimensional H\'enon map (\ref{3dHM1}) with fixed $B=-0.8$, $M_2=-1.05$ and varying $M_1$. When $-2.03 < M_1 < 2.172$, the fixed point $O$ is stable.
At $M_1\approx 2.172$ it undergoes a supercritical period-doubling bifurcation and  a stable period-2 cycle $(p_1,p_2)$ becomes attractor, Fig.~\ref{Lor_per2_bif}a. At $M_1 \approx 2.223$, the cycle $(p_1,p_2)$ undergoes a period doubling bifurcation: it becomes saddle and a stable period-4 cycle  $\tilde s = (s_1, s_2, s_3, s_4)$ emerges in its neighborhood, Fig.~10b.   A transition between Figs.~(b) and (c) seems insignificant, but in fact it is very important. Indeed, one can see that several events happened here: 1) a period-4 stable invariant curve $\tilde S = (S_1,S_2,S_3,S_4)$ are born from a stable period-4 cycle $\tilde s$; 2) separatrices of points $p_1$ and $p_2$ are reconstructed (e.g. the right separatrix of $p_1$ goes to the point $s_3$ at $M_1=2.27$, while, at $M_1=2.28$, it goes to the circle $S_1$); 3) the most important thing is that, at some $M_1$ between 2.27 and 2.28, a homoclinic butterfly of the cycle $(p_1,p_2)$ has appeared; 4) a saddle period-4 invariant curve $\tilde C = (C_1,C_2,C_3,C_4)$ is formed from this butterfly, since the saddle value of the cycle $(p_1,p_2)$ is greater than 1.
Almost immediately  the period-2 Lorenz attractor appear, see Fig.~\ref{Lor_per2_bif}d where it is seen that this attractor coexists with the stable period-4 invarinant curve $\tilde S$. When the curves
$\tilde S$ and $\tilde C$ merge and disappear, the period-2 Lorenz attractor becomes the unique attracting invariant set, see Figs.~\ref{Lor_per2_bif}e and \ref{Lor_per2_bif}(f), where such attractors are shown, without lacunae and with a lacuna, respectively.

Perhaps this will be very useful for the reader if we comment in more detail on some important  features of the structure of  the period-2 Lorenz attractor.

\begin{figure}
\centerline{\epsfig{file=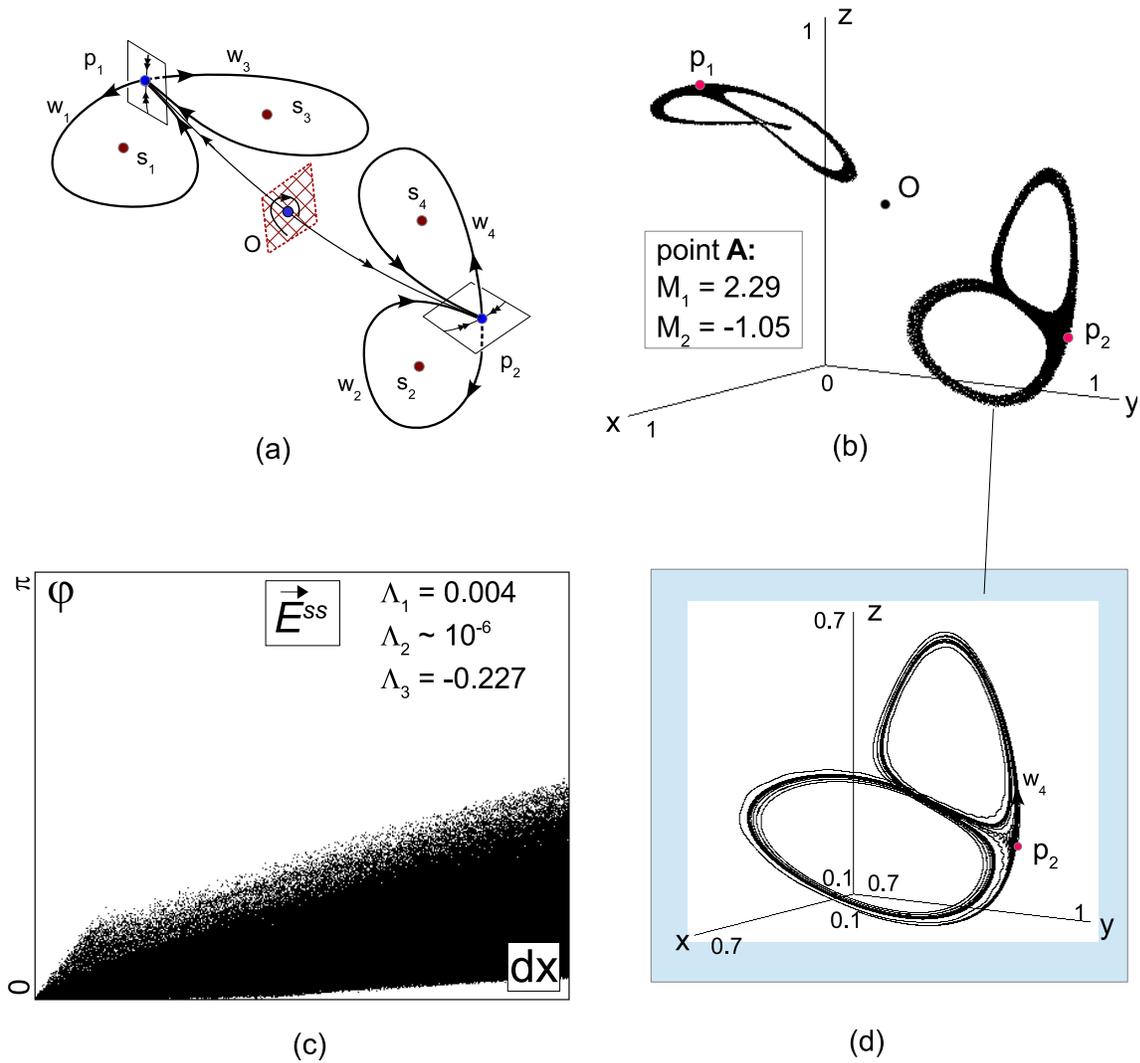, width=16cm
}}
\vspace{-1cm}
\caption{{\footnotesize Towards a structure of the period-2 Lorenz-like attractor: (a)  a skeleton scheme for the unstable separatrices of points $p_1$, $p_2$ and $O$;  (b) phase portrait for the attractor in map (\ref{3dHM1}) with $B=-0.8$ $M1 = -1.05$ and $M_2 = 2.29$; (c) its the LMP-graph, and
(d) numerics for its separatrix $w_4$.
Note also that, because of the period-2 Lorenz-like attractor is organized near a period-2 cycle, the LMP-graph in (c) was constructed by taking each 4th iteration of the map in order to check its pseudohyperbolicity in a proper way.
}}
\label{Lor2_skel_p_phyp}
\end{figure}

So, in Fig.~\ref{Lor2_skel_p_phyp}a, a skeleton scheme for the period-2 Lorenz-like attractor is shown that is similar to schemes from Figs.~\ref{DiscrLor}b and \ref{Discr-Lor_nor}b for the discrete orientable and nonorientable Lorenz attractors. The scheme in Fig.~\ref{Lor2_skel_p_phyp}a reflects main geometric properties of the attractor. So, the unstable manifold of the cycle $(p_1,p_2)$ is divided by the points $p_1$ and $p_2$ into 4 connected components, separatrices $w_1,w_2,w_3,w_4$,
such that $w_2 = T(w_1), w_3 = T(w_2), w_4 = T(w_3), w_1 = T(w_4)$. The two-dimensional stable manifold of the cycle $(p_1,p_2)$ consists also of two connected components, $W^s(p_1)$ and $W^s(p_2)$ that contain the points $p_1$ and $p_2$, respectively.  In the case under consideration, all separatrices $w_1,...,w_4$ do not intersect with $W^s(O)$ and, besides, the separatrices $w_1$ and $w_3$ intersect with $W^s(p_1)$ and do not intersect with $W^s(p_2)$, and the same holds symmetrically for $w_2$ and $w_4$. Importantly that the point $O$ and its unstable manifold do not belong to the attractor, and its stable manifold forms a natural boundary between two components of the period-2  Lorenz-like attractor (or a boundary between two discrete Lorenz attractor for the map $T^2$).

In Fig.~\ref{Lor2_skel_p_phyp}(b) and (d) there are shown pictures for the period-2 Lorenz-like attractor in map (\ref{3dHM1}) with $B=-0.8$, $M_1 = -1.05$ and $M_2 = 2.29$: the phase portrait of this attractor in Fig.~\ref{Lor2_skel_p_phyp}b and the separatrix $w_4$ in Fig.~\ref{Lor2_skel_p_phyp}d (the behavior of other separatrices, $w_1,w_2$ and $w_3$, is symmetric due to the points $p_1$ and $p_2$ compose a period-2 cycle and the unstable multiplier of this cycle is negative). It is important also that period-2 Lorenz-like attractors are often genuine pseudohyperbolic attractors. It concerns, in particular, to the attractors from Figs.~\ref{Lor_per2_bif}(d)--(f). For example, one can see from  the LMP-graph of Fig.~\ref{Lor2_skel_p_phyp}c constructed for the attractor of Fig.~\ref{Lor2_skel_p_phyp}b (it is also shown in Fig.~\ref{Lor_per2_bif}e) that this attractor is certainly pseudohyperbolic.


\subsection{On crises of period-2 Lorenz-like attractors.}

When change of parameters one can observe various types of crises with the period-2 Lorenz-like attractors. Usually, the same as the discrete Lorenz attractors \cite{GGS12,GGOT13,GGKT14}, they can transform into  closed invariant curves or strange attractors of torus-chaos type etc.
However, there are also own interesting features associated with the formation of attractor of new types. Below we discuss two of them.

\begin{figure}[ht]
\centerline{\epsfig{file=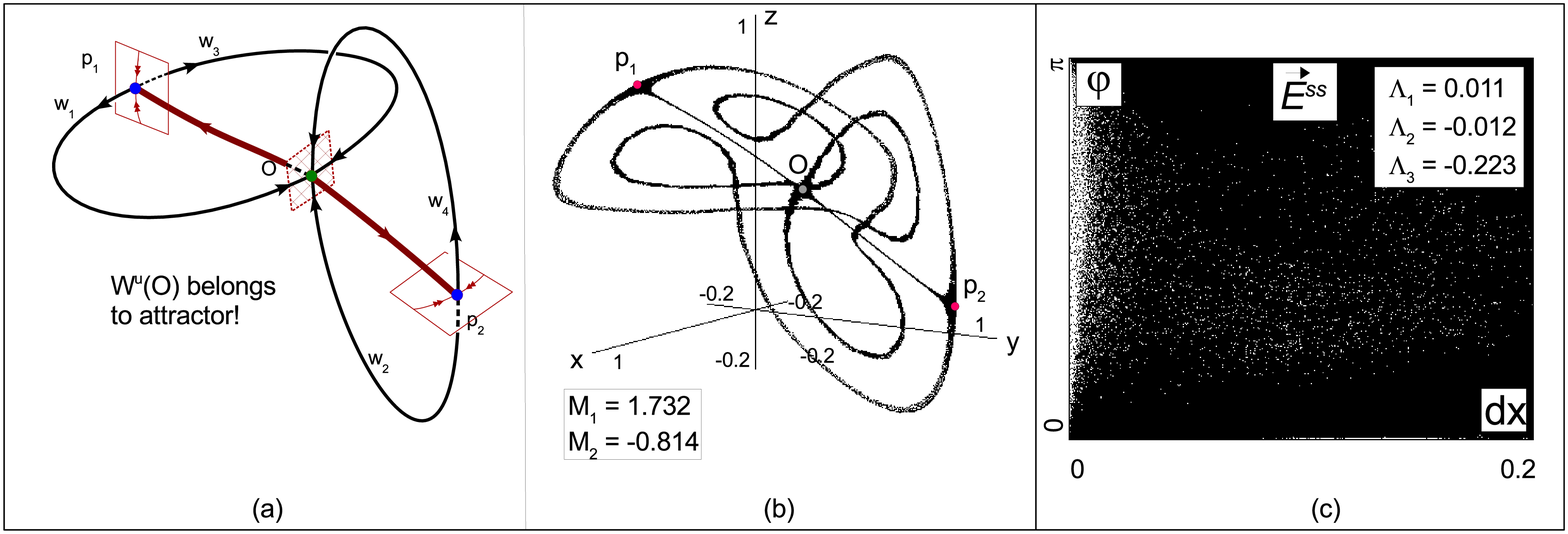, width=16cm
}}
\caption{{\footnotesize  Towards a homoclinic attractor containing the fixed point $O$ (that ia a saddle-focus near the 1:4 resonance) and the saddle cycle $(p_1,p_2)$ of period 2: (a) a skeleton scheme for such attractors; (b) an example of such attractor in map (\ref{3dHM1}) with $B=-0.8$, where $O=(x,x,x)$ with $x\approx 0.55$; (c) its LMP-graph. }}
\label{attr_skp_hom}
\end{figure}

When the unstable manifold of the cycle $(p_1,p_2)$ begin to intersect with $W^s(O)$, it follows that new intersections of $w_1$ and $w_3$ with $W^s(p_2)$ as well as $w_2$ and $w_4$ with $W^s(p_1)$ also appear. This means that a discrete homoclinic  attractor can arise containing the fixed point $O$ and the saddle cycle $(p_1,p_2)$ of period 2.   A skeleton scheme for such attractor  is shown in Fig.~\ref{attr_skp_hom}a. An example of such discrete attractor
was found in map (\ref{3dHM1}) with $B=-0.8$ at $M_1 = 1.732, M_2 = -0.814$, see Fig.~\ref{attr_skp_hom}b. This attractor contains the points $O$, $p_1$ and $p_2$ and entirely their unstable invariant manifolds that have structures like ``wings'' whose skeleton is formed pairwise from separatices  $w_1,w_3$ and $w_2,w_4$, respectively, see Fig.~\ref{attr_skp_hom}a.
Note that this attractor in the case of Fig.~\ref{attr_skp_hom}b contains
%
the fixed point $O$
with multipliers $\gamma<-1$ and $\lambda_{1,2} = \rho e^{\pm i\varphi}$, where $\rho\approx 0.81$ and $\varphi$ is very close to $\pi/2$ ($\cos\varphi \approx 0.03$), i.e. we have a situation near the
strong resonance 1:4.
Thus, we also touch upon the problem of the structure of arising homoclinic attractors when passing near strong resonances, which was formulated in \cite{GGS12}.\footnote{In particular, some interesting types of discrete spiral attractors were found in three-dimensional generalized H\'enon maps near strong resonances: the so-called discrete super-spiral attractor (containing period-4 saddle-foci of both types (2,1) and (1,2)) was found in \cite{GGS12} nearly the strong resonance 1:4, and a discrete ``triangle'' spiral homoclinic attractor was found in \cite{GGKT14} near the strong resonance 1:3.} The attractor from Fig.~\ref{attr_skp_hom}b is certainly a quasiattractor, since it contains the fixed point  $O$ that is a saddle-focus (2,1)
%
%
which, besides, has the saddle value less than 1. Indirectly, this is confirmed by its LMP-graph, see Fig.~\ref{attr_skp_hom}c.

\begin{figure}[ht]
\centerline{\epsfig{file=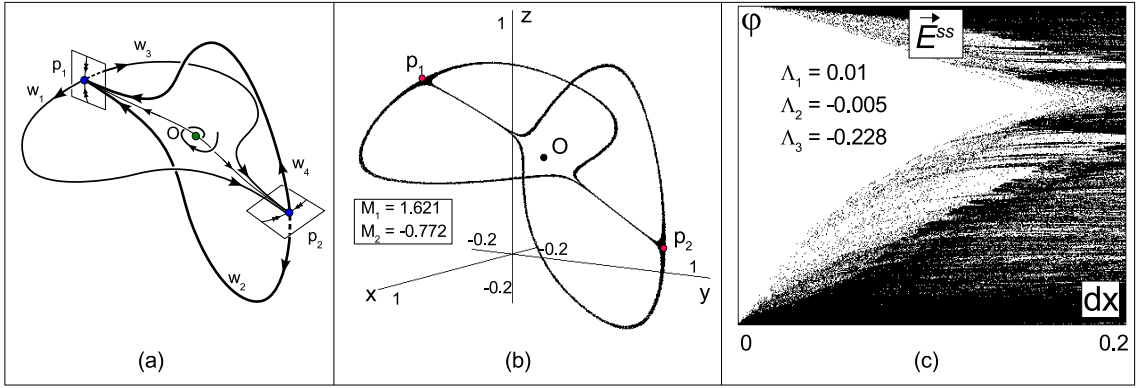, width=16cm
}}
\vspace{-1cm}
\caption{{\footnotesize Towards a period-2 heteroclinic attractor containing the the saddle cycle $(p_1,p_2)$: (a) a skeleton scheme for such attractors; (b) an example of such attractor in map (\ref{3dHM1}) with $B=-0.8$, where $p_1 =(a,b,a), p_2 = (b,a,b)$ and $a\approx 0.85, b\approx 0.12$; (c) its LMP-graph.}}
\label{attr_skp_het}
\end{figure}

One can imagine such situation, when, with varying parameters, the unstable separatrices of the cycle $(p_1,p_2)$ do not longer intersect with the stable manifold $W^s(O)$ of the fixed point $O$ and,
thus, a heteroclinic period-2 attractor is created.
It contains the cycle $(p_1,p_2)$ of period 2 but does not capture the fixed point $O$. A skeleton scheme of such attractor is shown in Fig.~\ref{attr_skp_het}a, it illustrates the main feature of attractor related to the fact that all stable and unstable invariant manifolds of points $p_1$ and $p_2$ are mutually intersect. An example of such attractor for map (\ref{3dHM1}) with $B=-0.8$ is shown in Fig.~\ref{attr_skp_het}b. However, we  note that the first example of such attractor was found yet in \cite{GOST05}, see Fig.~\ref{fig:attrmap1}d. But now we understand better both its structure and how it can appear at reconstructions with the period-2 Lorenz attractors. It is also important that such period-2 heteroclinic attractor can be genuine, pseudohyperbolic.
Indirectly, this seems to confirm its LMP-graph from Fig.~\ref{attr_skp_het}c, but here we can not say with great certainty that the attractor is pseudohyperbolic indeed. Although we construct this graph using every 4th iteration of the map,  a certain closeness of the graph to the point $(0,\pi)$ do not allow to predict that a tendency of approaching graph to this point will not terminate at increasing time of calculations. So, we consider the question of pseudohyperbolicity of the attractor from Fig.~\ref{attr_skp_het}c as a still open problem.


\begin{figure}[ht]
\centerline{\epsfig{file=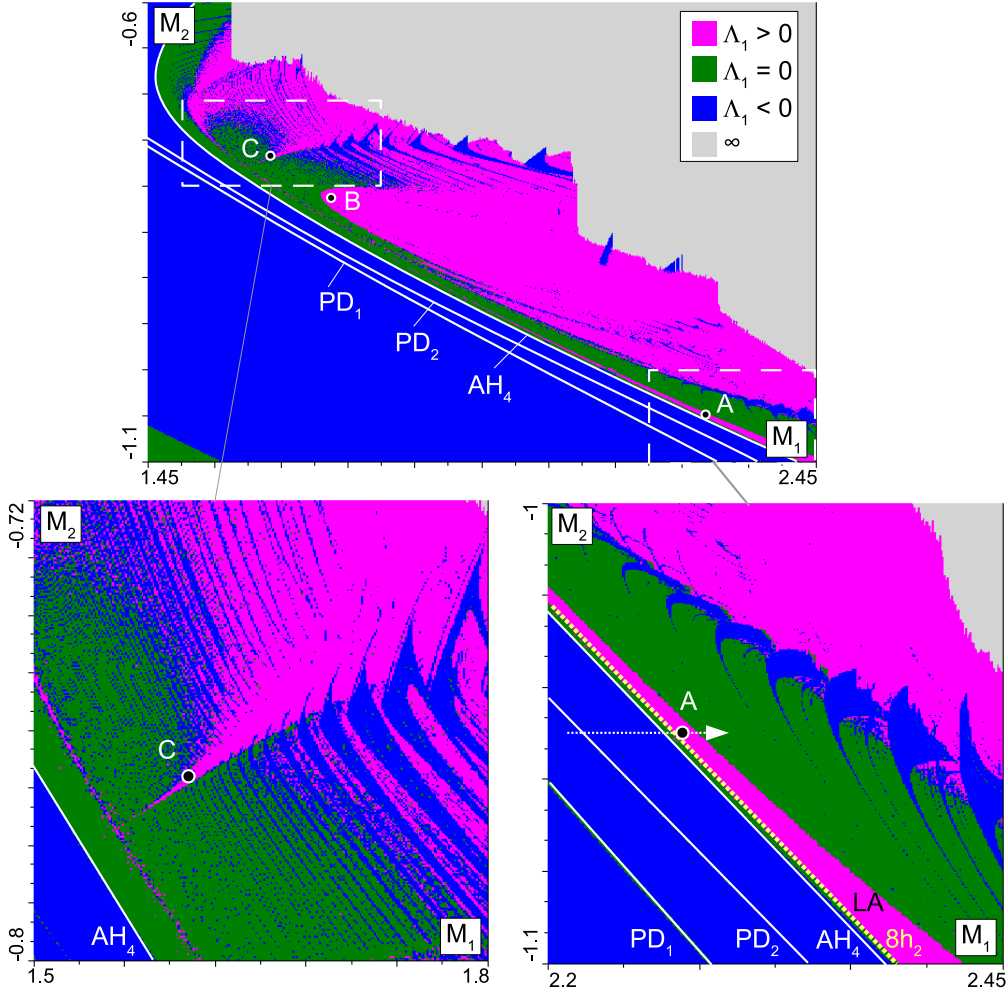, width=14cm
}}
\vspace{-1cm}
\caption{{\footnotesize Chart of Lyapunov exponents for map (\ref{3dHM1}) with $B = - 0.8$ on the rectangle $1.45 < M_1 < 2.45; -1.1 < M_2 < -0.6$ of the $(M_1,M_2)$-parameter plane. 
Some bifurcation curves are also shown here: the curves $PD_1$ and $PD_2$ of the first and second period doubling with the fixed point $O$; the curve $AH_4$ of a supercritical Andronov-Hopf bifurcation of the stable period-4 cycle $(s_1,s_2,s_3,s_4)$; and the curve (a think zone indeed) $8h_2$ corresponding to the appearance of a homoclinic butterfly at the cycle $(p_1,p_2)$. }}
\label{B_m0p8_LE}
\end{figure}

In order to trace better interrelations between the discrete attractors of these three types, the period-2 Lorenz attractor, the homoclinic attractor with period-2 cycle and the period-2 heteroclinic attractor,  we have constructed the so-called chart of maximal Lyapunov exponents, see Fig.~\ref{B_m0p8_LE}, for map (\ref{3dHM1}) with $B = - 0.8$. This chart is a specific diagram, on the $(M_1,M_2)$-parameter plane, showing various types of stable regimes.
%
The blue and green domains in the chart relate to stable periodic orbits and stable closed invariant curves, respectively. Strange attractors exist for the values of parameters corresponding to the rose domains. In the latter domains we mark three points A, B and C related to the attractors presented in Figs.~\ref{Lor2_skel_p_phyp},~\ref{attr_skp_hom} and~\ref{attr_skp_het}, respectively.

Even a quick glance at the diagram of Fig.~\ref{B_m0p8_LE} is enough to understand that the transitions between different attractors corresponding to points A, B and C can not be too simple. In any case, they are associated not only with the described above reconstructions of unstable manifolds but also, mostly, with a destruction of attractors themselves, their transformations into invariant curves, inverse rearrangements of these curves into new attractors, etc. Understanding the corresponding bifurcation mechanisms even in the particular case of map (\ref{3dHM1}) looks potentially as a new rather interesting task. In the present paper we considered, in fact, only one such mechanism -- the scenario of appearance of the period-2 Lorenz attractor due to the option {\bf[sc3]} that is observed in the way marked by the turquoise arrow in Fig.~\ref{B_m0p8_LE}.

\section{Conclusion}  \label{sec:concl}

Discrete homoclinic attractors of multidimensional maps compose a very interesting class of new strange attractors that can be genuine, pseudohyperbolic, attractors. The theory of such attractors has not been created yet, although some important results in this direction were obtained, see e.g. \cite{GOST05,GGS12,GGOT13,GGKT14}. The greatest advances in this theory have been achieved at the study of discrete Lorenz-like attractors. In the present paper we tried to give some overview of main results obtained in this direction.
However, as one can see, even here there are still many unsolved problems. For example, can there exist in  three-dimensional H\'enon maps discrete Lorenz attractors containing saddles with positive multipliers? Can discrete Lorenz-like attractors of nonorientable maps in three-dimensional generalized H\'enon maps and be pseudohyperbolic? For the latter, it is also not yet known whether they can be born under local bifurcations of fixed points with a triplet $(+1,+1,-1)$ of multipliers, etc.

On the other hand, as we know, the zoo of discrete homoclinic attractors even in three-dimensional maps now contains several specimens, interesting  from different points of view. These are, for example, a series of figure-8 attractors, including double and super figure-8 ones, as well as figure-8 spiral attractors, see \cite{GG16}). Concerning discrete Lorenz-like attractors, some their types, predicted by their primary homoclinic structures, do not found yet. These are, e.g. mentioned above ``classical'' attractors (containing a saddle (2,1) point with all positive multipliers) or nonorientable ``double Lorenz'' attractors (containing a saddle (2,1) point with multipliers $\gamma>1$,
$-1<\lambda_2<0<\lambda_1<1$).  Very interesting class consists of the so-called discrete Shilnikov attractors, i.e. homoclinic attractors containing a saddle-focus fixed point of type (1,2). As such attractors contain the closure of two-dimensional unstable manifold of this point, then they can give quite realistic and simple examples of hyper-chaotic attractors \cite{GSKK19,Stan19}. Some results related to the the theory of discrete Shilnikov attractors were obtained e.g. in \cite{GGS12,GGKT14,GG16,AGK17}. However, there is still much that can be done.

Another interesting observation can be made in continuation of the topics developed in Section~\ref{sec:othertype} of this paper. Namely, in a generalization of the options {\bf[sc1]}, {\bf[sc2]} and {\bf[sc3]}, we consider the following\\

\begin{itemize}

\item[{\bf[sc$_n$]}]
the fixed point $O$
undergoes $n$ supercritical period doubling  bifurcations and this sequence is terminated by the Andronov-Hopf bifurcation, i.e. the stable cycle of period $2^n$
loses stability according to the options [s1] or [s2] for $T^{2n}$.\\

\end{itemize}

As a result of this option, a period-$2^n$ Lorenz-like attractor can theoretically appear. Next, components of this attractor can merge pairwise forming a $2^{n-1}$-component attractor containing period-$2^n$ points and organizing them (after period doubling) period $2^{n-1}$-points, etc.

This option implies that, in principle, such a sequence of attractor crises can be observed:    a period-$2^n$ Lorenz-like attractor appears $\Rightarrow$  $2^n$ components of this attractor merge in pairs forming $2^{n-1}$-component attractor $\Rightarrow$ ... $\Rightarrow$ a homoclinic (of Lorenz shape or not) attractor containing the fixed point $O$ appears. This sequence of attractor crises can be viewed as a ``chain of doublings of discrete Lorenz-shape attractors''.
In this paper, we only slightly touched this topic on the example of such chain with n = 1.
Of course, the question on such chains with larger
values of $n$ seems rather interesting, and we plan to consider this problem in more detail in the nearest future.

\section*{Acknowledgements}
This paper was carried out in the framework of RSF grants 18-71-00127 (Sec. 2 and 3) and 19-71-10048 (Sec. 4 and 5). The authors are partially supported (numerical results of Section~\ref{LMPmethod}) by Laboratory of Dynamical Systems and Applications NRU HSE, of the Ministry of science and higher education of the RF grant ag. No. 075-15-2019-1931.


\end{document}